\documentclass[11pt]{article}

\usepackage{amscd,amsmath, amssymb, fancyhdr, epsfig, color, tocloft, mathtools}
\usepackage{dutchcal, extarrows}

\usepackage[backref=page]{hyperref}
\renewcommand*{\backrefalt}[4]{%
	\ifcase #1 (Not cited.)%
	\or        (Cited on page~#2.)%
	\else      (Cited on pages~#2.)%
	\fi}

\hypersetup{
	colorlinks   = true,
	citecolor    = magenta,
	linkcolor    = blue,
	urlcolor     = magenta	
}

\newcommand{\version}{version 2.0,\ \ Nov. 14, 2024}
\makeatletter
\def\x@arrow{\DOTSB\Relbar}
\def\xlongequalsignfill@{\arrowfill@\x@arrow\Relbar\x@arrow}
\providecommand{\xlongequal}[2][]{%
	\ext@arrow 0099\xlongequalsignfill@{#1}{#2}}
\def\xlongrightarrowfill@{\arrowfill@\relbar\relbar\longrightarrow}
\makeatother

\numberwithin{equation}{section}

\def\eqref#1{(\ref{#1})}

\newcommand{\g}{{\mathfrak g}}

\newcommand{\Z}{{\mathbb Z}}
\newcommand{\C}{{\mathbb C}}
\newcommand{\R}{{\mathbb R}}
\newcommand{\Q}{{\mathbb Q}}

\newcommand{\6}{\partial}
\def\1{\sqrt{-1}\:}

\newcommand{\cntrct}                
{\hspace{2pt}\raisebox{1pt}{\text{$\lrcorner$}}\hspace{2pt}}
\newcommand{\arrow}{{\:\longrightarrow\:}}

\newcommand{\calo}{{\cal O}}

\renewcommand{\bar}{\overline}
\renewcommand{\phi}{\varphi}
\renewcommand{\epsilon}{\varepsilon}
\renewcommand{\geq}{\geqslant}

\newcommand{\im}{\operatorname{im}}

\newcommand{\Vol}{\operatorname{Vol}}

\newcommand{\Av}{\operatorname{Av}}

\renewcommand{\Re}{\operatorname{Re}}
\renewcommand{\Im}{\operatorname{Im}}

\newcounter{Mycounter}[section]
\newcounter{lemma}[section]
\newcounter{claim}[section]
\renewcommand{\theclaim}{{Claim \thesection.\arabic{claim}}}
\newcommand{\claim}{%
	\setcounter{claim}{\value{Mycounter}}
	\refstepcounter{claim}
\stepcounter{Mycounter}
	{\noindent \bf \theclaim:\ }}

\newcounter{sublemma}[section]
\newcounter{corollary}[section]
\renewcommand{\thecorollary}{{Corollary \thesection.\arabic{corollary}}}
\newcommand{\corollary}{%
	\setcounter{corollary}{\value{Mycounter}}
	\refstepcounter{corollary}
	\stepcounter{Mycounter}
	{\noindent \bf \thecorollary:\ }}

\newcounter{theorem}[section]
\renewcommand{\thetheorem}{{Theorem \thesection.\arabic{theorem}}}
\newcommand{\theorem}{%
	\setcounter{theorem}{\value{Mycounter}}
	\refstepcounter{theorem}
	\stepcounter{Mycounter}
	{\noindent \bf \thetheorem:\ }}

\newcounter{conjecture}[section]
\renewcommand{\theconjecture}{{Conjecture \thesection.\arabic{conjecture}}}
\newcommand{\conjecture}{%
	\setcounter{conjecture}{\value{Mycounter}}
	\refstepcounter{conjecture}
	\stepcounter{Mycounter}
	{\noindent \bf \theconjecture:\ }}

\newcounter{proposition}[section]
\renewcommand{\theproposition} {{Proposition \thesection.\arabic{proposition}}}
\newcommand{\proposition}{%
	\setcounter{proposition}{\value{Mycounter}}
	\refstepcounter{proposition}
	\stepcounter{Mycounter}
	{\noindent \bf \theproposition:\ }}

\newcounter{definition}[section]
\renewcommand{\thedefinition} {{Definition~\thesection.\arabic{definition}}}
\newcommand{\definition}{%
	\setcounter{definition}{\value{Mycounter}}
	\refstepcounter{definition}
	\stepcounter{Mycounter}
	{\noindent \bf \thedefinition:\ }}

\newcounter{example}[section]
\renewcommand{\theexample}{{Example \thesection.\arabic{example}}}
\newcommand{\example}{%
	\setcounter{example}{\value{Mycounter}}
	\refstepcounter{example}
	\stepcounter{Mycounter}
	{\noindent \bf \theexample:\ }}

\newcounter{remark}[section]
\renewcommand{\theremark}{{Remark \thesection.\arabic{remark}}}
\newcommand{\remark}{%
	\setcounter{remark}{\value{Mycounter}}
	\refstepcounter{remark}
	\stepcounter{Mycounter}
	{\noindent \bf \theremark:\ }}

\newcounter{problem}[section]
\newcounter{question}[section]
\makeatletter

\@addtoreset{equation}{section}
\@addtoreset{footnote}{section}
\makeatother

\def\blacksquare{\hbox{\vrule width 5pt height 5pt depth 0pt}}
\def\endproof{\blacksquare}

\newcommand{\proof}{{\bf Proof: \ }}
\newcommand{\pstep}{{\bf Proof. Step 1: \ }}

\begin{document}
	
\begin{center}
{\Large\bf  The Lee--Gauduchon cone on complex manifolds}\\[5mm]

		{\large
			Liviu Ornea\footnote{Liviu Ornea is  partially supported by the PNRR-III-C9-2023-I8 grant CF 149/31.07.2023 Conformal
				Aspects of Geometry and Dynamics.},  
			Misha Verbitsky\footnote{Misha Verbitsky is partially supported by
				FAPERJ 	SEI-260003/000410/2023 and CNPq - Process 310952/2021-2.\\[1mm]
			}\\[4mm]
			
		}
		
	\end{center}

	\hfill
	
{\small
\hspace{0.15\linewidth}
\begin{minipage}[t]{0.7\linewidth}
{\bf Abstract} \\
Let $M$ be a compact complex $n$-manifold. A Gauduchon
metric is a Hermitian metric whose fundamental 2-form
$\omega$ satisfies the equation
$dd^c(\omega^{n-1})=0$. Paul Gauduchon has proven that any
Hermitian metric is conformally equivalent to a Gauduchon
metric, which is unique (up to a constant multiplier) in
its conformal class. Then $d^c(\omega^{n-1})$ is a closed
$(2n-1)$-form; the set of cohomology classes of all such
forms, called the Lee-Gauduchon cone, is a convex
cone, superficially similar to the K\"ahler cone.  
We prove that the Lee-Gauduchon
cone is a bimeromorphic invariant, and compute it 
for several classes of non-K\"ahler manifolds.
\end{minipage} 
	}

\tableofcontents
	

\section{Introduction}	


This paper was inspired by our work \cite{_OV_balanced_} on 
the following conjecture.

\hfill

\conjecture
A compact locally conformally K\"ahler manifold
never admits a balanced Hermitian metric.

\hfill

Both of these notions are defined properly
in Section \ref{_examples_Section}.

\hfill

It turns out that this question
is essentially a question about the shape
of the Gauduchon cone, the cone which plays
in non-K\"ahler complex geometry the same role
as the K\"ahler cone plays in the conventional (K\"ahler)
complex geometry. 

A Hermitian metric $\omega$ on a complex $n$-manifold
is called {\bf Gauduchon} if $dd^c \omega^{n-1}=0$
(see Subection \ref{_Gauduchon_cone_Subsection_}
for a proper introduction to the
geometry of Gauduchon forms).
The {\bf Gauduchon form} $\omega^{n-1}$
can be used to reconstruct $\omega$ in an
unambigous way. The set of all Gauduchon forms
is clearly a convex cone in the infinite-dimensional
space of all $dd^c$-closed $(n-1,n-1)$-forms.

To make sense of the cohomological properties
 of the Gauduchon forms, one uses the Aeppli
cohomology, defined as $\frac{\ker dd^c}{\im d + \im d^c}$
(Subsection \ref{_Aeppli_Subsection_}).
The Gauduchon metrics define a convex cone,
called {\bf the Gauduchon cone},
in the Aeppli cohomology group $H^{n-1,n-1}_{AE}(M)$.

For a general complex manifold,
the Aeppli cohomology is hard to compute;
this group is not a topological invariant,
and not constant in holomorphic families
of complex manifolds. In this paper we 
introduce a more accessible counterpart
of the Gauduchon cone, called the 
Lee-Gauduchon cone, which lives
in the de Rham cohomology group $H^{2n-1}(M,\R)$.
We introduce this object and give its
proper definition in Subsection \ref{_LG_cone_Subsection_}.
From the definition of a Gauduchon form,
it is clear that the form $d^c(\omega^{n-1})$ is 
closed. The Lee-Gauduchon cone is the
set of all cohomology classes of this form.

In this paper we show that the Lee-Gauduchon
cone is a bimeromorphic invariant (\ref{_LG_bimeromo_Theorem_}) and compute
it for several important classes of complex manifolds.
We also characterize the Lee-Gauduchon cone
in terms of positive currents (\ref{_LG_dual_to_pos_Theorem_}).

\section{Bott--Chern and Aeppli cohomologies}

We present briefly the two main cohomology theories which
are most useful  on non-K\"ahler manifolds.
We refer to \cite{_Angella:book_} for more details and reference.

\subsection{Bott--Chern cohomology}

\definition
Let $M$ be a complex manifold, and $H^{p,q}_{BC}(M)$
the space of closed $(p,q)$-forms modulo $dd^c(\Lambda^{p-1,q-1} M)$.
Then $H^{p,q}_{BC}(M)$ is called {\bf the Bott--Chern
	cohomology} of $M$. 
	
\hfill

\remark
A $(p, q)$-form $\eta$ is closed if and only if
$\6\eta=\bar\6\eta=0$. 
Since $2\1 \6\bar\6= dd^c$, an equivalent definition of the 
	Bott--Chern cohomology is  
	$H^*_{BC}(M):=\frac{\ker \6 \cap \ker \bar\6}{\im \6\bar\6}$.
	
\hfill

\remark
While there exist no morphisms between the de Rham and
the Dolbeault cohomology, there are natural and functorial maps from the
Bott--Chern cohomology to both the Dolbeault cohomology
$H^*(\Lambda^{*,*} M, \bar\6)$ and  the de Rham
cohomology. However, there exists no multiplicative structure on 
the Bott--Chern cohomology.

\hfill


Denote by $H^{1,0}_{d} (M)$ the space of closed
holomorphic forms on $M$. Since any exact holomorphic 1-form
on a compact complex manifold vanishes, $H^{1,0}_{d} (M)$ can be considered as a
subspace in de Rham (or Dolbeault) cohomology.
Clearly, $H^{1,0}_{d} (M)\oplus \overline{H^{1,0}_{d} (M)}$
is the space of all forms $\alpha\in \Lambda^1(M, \C)$ such that
$d\alpha=d^c\alpha=0$; this space is the complexification 
of $\ker d\oplus \ker d^c$.

\hfill

%
%
%
%
%
%

\proposition
The sequence 
\begin{equation}\label{_3rd_exact_sequence_in_BC_} 
0\arrow \Re(H^{1,0}_{d} (M)) 
\arrow H^1(M, \R) 
\stackrel{d^c}{ \arrow }
	H^{1,1}_{BC}(M,\R)
\end{equation}
is exact.

\proof
This sequence is clearly exact in the first term:
an exact holomorphic form is a differential of a global holomorphic
function on $M$, and all such functions are constant. To prove that
it is exact in the second term, let $x$ be a closed 1-form, and $[x]\in H^1(M)$
its cohomology class. The cohomology class of $d^c x$ vanishes in $H^{1,1}_{BC}(M,\R)$
if and only if $d^c x= dd^c f$, for some function $f\in C^\infty M$.
However, $d^c x= dd^c f$ means that $x+ df$ is $d$-closed and $d^c$-closed,
hence $[x]$ belongs to the image of $H^{1,0}_{d} (M)\oplus \overline{H^{1,0}_{d} (M)}$.
\endproof

%
%
%

\subsection{Aeppli cohomology}
\label{_Aeppli_Subsection_}

\definition 
Let $M$ be a complex manifold, and $H^{p,q}_{AE}(M)$
the space of $dd^c$-closed $(p,q)$-forms modulo
$\6(\Lambda^{p-1,q} M)+ \bar\6(\Lambda^{p,q-1} M)$.
Then $H^{p,q}_{AE}(M)$ is called {\bf the Aeppli 
	cohomology} of $M$. 

\hfill

\remark
The spaces $H^{p,p}_{AE}(M)$ and $H^{p,p}_{BC}(M)$
are preserved by the complex conjugation. We denote
the subspaces fixed by the complex conjugation by
$H^{p,p}_{AE}(M, \R)\subset H^{p,p}_{AE}(M)$ and
$H^{p,p}_{BC}(M,\R)\subset H^{p,p}_{BC}(M)$.	

\hfill

\theorem (A. Aeppli, \cite{_Aepli:sequences_})  
Let $M$ be a compact complex $n$-manifold. Then
{the Aeppli cohomology is finite-dimensional.} Moreover, 
the natural pairing 
\begin{equation*}
	\begin{split}
		 H^{p,q}_{BC}(M) \times H^{n-p,n-q}_{AE}(M)&\arrow
H^{2n}(M)=\C,\\
 (x,y) & \mapsto\int_M x\wedge y,
\end{split}
\end{equation*}
{is non-degenerate and identifies $H^{p,q}_{BC}(M)$
	with the dual $(H^{n-p,n-q}_{AE}(M))^*$.}

\subsection{The Lee-Gauduchon space 
$W\subset H^{2n-1}(M)$}
\label{_W:subsection_}

We define a subspace in the de Rham cohomology which will play a central role in the sequel.
The motivation for its name will be apparent after we define the Lee-Gauduchon cone. 

\hfill

\definition Let
$W\subset H^{2n-1}(M, \R)$ be the space of all cohomology classes $\alpha$
such that $\int_M \alpha \wedge \rho=0$
for all closed holomorphic forms $\rho\in
\Lambda^{1,0}(M, \R)$. We call $W$ {\bf the Lee-Gauduchon space}.

\hfill

\remark
Since $\int_M d^c u \wedge \rho= 
\int_M \omega^{n-1} \wedge d^c \rho= 0$,
the image of the natural map 
$d^c:\;H^{n-1,n-1}_{AE}(M)\arrow H^{2n-1}(M)$ belongs to $W$.

\hfill

\proposition \label{_W_in_terms_of_AE:Proposition_}
Let $M$ be a compact complex $n$-manifold.
Then  
\[ W=d^c(H^{n-1,n-1}_{AE}(M, \R)).
\]

\proof 
Consider the exact sequence \eqref{_3rd_exact_sequence_in_BC_} 
constructed above. Dualizing it, we obtain
\[
H^{n-1,n-1}_{AE}(M, \R) \stackrel{d^c}\arrow H^{2n-1}(M, \R)
\arrow H^{2n-1}(M, \R)/W \arrow 0,
\]
because $W$ is the annihilator of $H^{1,0}_{d}\oplus \overline{H^{1,0}_{d}
	(M)}\subset H^1(M, \R)$.
\endproof

\subsection{The Gauduchon cone of a complex manifold}
\label{_Gauduchon_cone_Subsection_}

\definition\label{_Gauduchon_metric_Definition_}
A Hermitian metric on a complex $n$-manifold is 
called {a \bf Gauduchon metric} if its fundamental form satisfies  $dd^c(\omega^{n-1})=0$. In this case, $\omega^{n-1}$ is called {\bf a Gauduchon form}.

\hfill

\remark Recall that on a compact complex manifold, a Gauduchon metric exists in any conformal class of Hermitian metrics, and it is unique up to constant multiplier, \cite{_Gauduchon_}.
	
\hfill

\definition
{\bf The Gauduchon cone} of a compact complex
$n$-manifold is the set of all classes
$\omega^{n-1}\in H^{n-1, n-1}_{AE}(M, \R)$ of all Gauduchon
forms.

\hfill 

\remark
Since any strictly positive $(n-1, n-1)$-form is the 
$(n-1)$-th power of a Hermitian form (\cite{_Michelson_}), 
the Gauduchon cone is open and convex in $H^{n-1, n-1}_{AE}(M, \R)$.

\subsection{Positive currents}

In this subsection we express the Gauduchon cone in terms of currents. We recall the definition of positive currents and refer the reader to \cite{_Demailly:CADG_} for an introduction on currents on complex manifolds. 

\hfill

%

\definition\label{_Positive_current:Definition_}
A $(1,1)$-current $\alpha$ is called {\bf positive}
if $\int_M \alpha \wedge \tau \geq 0$ for
any positive $(n-1,n-1)$-form $\tau$ with compact support.

\hfill 

\remark
The cone of positive $(1,1)$-currents is generated by
$-\1 \alpha u\wedge \bar u$, where $\alpha$ is a
positive generalized function (that is, a measure),
and $u$ a (1,0)-form.

\hfill

\definition
The {\bf  pseudo-effective cone}
$P\subset  H^{1,1}_{BC}(M)$
is the set of Bott--Chern classes of all
positive, closed (1,1)-currents.

\hfill

The next result gives a characterization of  the Gauduchon cone in terms of currents.

\hfill

\theorem  \label{_Lamari_pseff_Theorem_}
On a compact complex manifold, the Gauduchon cone 
is dual to the pseudo-effective cone.

\hfill

\proof \cite[Lemme 3.3]{_Lamari_}; see also \cite{_Popovici_Ugarte_}.
We include the proof of \ref{_Lamari_pseff_Theorem_}
in Section \ref{_exact_pseff_Subsection_}.
 \endproof

\section{Lee forms and Lee classes of a Gauduchon form}

\subsection{The Lee-Gauduchon cone}
\label{_LG_cone_Subsection_}

The following is but one of the equivalent definitions of the Lee form of a Hermitian metric. For another approach, see \cite{_Gauduchon_}.

\hfill
	
\definition
	{\bf The Lee form} of a Hermitian metric
	$\omega$ is $\frac 1 {n-1} *(d^c\omega^{n-1})$,
where $*$ denotes the Hodge star operator.
	
%
	
\hfill

\remark
	Clearly, 
	the Lee form is $d^*$-closed 
	{if and only if $\omega$ is Gauduchon.}
	
\hfill

\definition
	Let $\omega^{n-1}$ be a Gauduchon form. The corresponding 
	{\bf Lee-Gauduchon form} is
	$d^c (\omega^{n-1})$. This form is clearly closed;
	its cohomology class $[\omega^{n-1}] \in H^{2n-1}(M,\R)$ is called {\bf the Lee-Gauduchon class}.

\hfill

\remark
	As in Subsection \ref{_W:subsection_}, let $W\subset H^{2n-1}(M,\R)$
	be the space of all cohomology classes $\alpha$
	such that $\int_M \alpha \wedge \rho=0$
	for all closed holomorphic forms $\rho\in \Lambda^{1,0}(M)$.
	Since $\int_M d^c (\omega^{n-1}) \wedge \rho= 
	\int_M \omega^{n-1} \wedge d^c \rho= 0$,
	all Lee-Gauduchon classes belong to $W$, called {\bf the Lee-Gauduchon space}
(Subsection \ref{_W:subsection_}).

\hfill

\definition
	{\bf The Lee-Gauduchon cone} $LG(M)\subset W$ is 
	the set of all Lee-Gauduchon classes of Gauduchon forms.

\hfill

\claim
	The Lee-Gauduchon cone is a convex, open cone in the Lee-Gauduchon space $W$.\\
	
\proof
	The set of Lee-Gauduchon forms is open,
	because it is the image of the Gauduchon cone, which is open,
	and the projection from $dd^c$-closed forms to $W$ is
	surjective by \ref{_W_in_terms_of_AE:Proposition_}.
	\endproof

\hfill
	
\remark
The Lee-Gauduchon space is trivial on all compact complex manifolds for which the $dd^c$-lemma
		holds, for example, on all projective, Moishezon, 
	        K\"ahler or Fujiki class C manifolds.

\subsection{Exact pseudoeffective Bott--Chern classes}
\label{_exact_pseff_Subsection_}

In this subsection we express the Lee-Gauduchon cone in
terms  of exact pseudoeffective Bott--Chern classes.

We use the Hahn--Banach theorem 
(on the model of \cite{_Sullivan:foliated_complex_,_HL:intrinsic_}) 
in the following formulation:

\hfill

\theorem (\rm Hahn--Banach)\label{_Hahn_Banach:Theorem_} 
Let $V_1$ be a locally convex topological vector space,
$V\subset V_1$ a closed subspace, and $A\subset V_1$ an open,
convex subset, not intersecting $A$. Then there exists
	a continuous linear functional $\xi\in V^*_1$ vanishing on $V$ and
	positive on $A$.
	
\hfill

We start by proving Lamari's theorem
(\ref{_Lamari_pseff_Theorem_}).

\hfill

\definition
Let $\alpha \in \Lambda^{p,q}(M)$ be a $dd^c$-closed
form. We say that $\alpha$ is {\bf Aeppli exact}
if it is the $(p,q)$-part of an exact form.

\hfill

\theorem\label{_Lamari_proved_Theorem_}
On a compact complex manifold, the Gauduchon cone 
is dual to the pseudo-effective cone.

\hfill

\proof
Let $A$ be the set of all strictly 
positive $(n-1, n-1)$-forms, $u\in H^{n-1,n-1}_{AE}(M, \R)$,
and $V=u + V'$, where $V'$ is the space of Aeppli exact $dd^c$-closed 
$(n-1, n-1)$-forms. Then $V \cap A=\emptyset$
if and only if there exists a functional on $\Lambda^{n-1,n-1}(M)$
(that is, a (1,1)-current) $\xi$ such that $\langle \xi, A \rangle > 0$
and $\langle \xi, V \rangle = 0$. Since $V$ is an affine space and $V'$
its linearization, $\langle \xi, V \rangle = 0$ implies $\langle \xi, V' \rangle = 0$.

The condition $\langle \xi, A \rangle > 0$
means precisely that $\xi$ is a non-zero positive current. 
The condition $\langle \xi, V' \rangle = 0$
is equivalent to $\int \xi \wedge dw=0$
for all $(2n-3)$-forms $w$, because $V'$ is the space of
$(n-1,n-1)$-parts of exact forms, and $\xi$ is a (1,1)-current.
However, $\int \xi \wedge dw=0$ for all $w$ is equivalent
to $\xi$ being closed. Then, a class $u\in H^{n-1,n-1}_{AE}(M, \R)$
belongs to the Gauduchon cone if and only if $\int_M \xi \wedge u >0$
for all positive, closed (1,1)-currents $\xi$.
\endproof

\hfill 

%
%

We may now prove:

\hfill

\theorem \label{_LG_dual_to_pos_Theorem_}
Let $M$ be a compact complex manifold,
and ${\cal C}\subset H^1(M)$ the set of classes
$\rho\in H^1(M,\R)$ such that $d^c(\rho)\in H^{1,1}_{BC}(M)$
is pseudo-effective. Then the Lee-Gauduchon
	cone $LG(M)\subset W$ is 
	the dual cone to ${\cal C}$. In other words,
 $\alpha \in  LG(M)$ if and only if $\int_M \alpha \wedge \rho>0$
	for any closed 1-current $\rho$ such that $d^c\rho$ is positive.

\hfill
	
\pstep
If $\alpha \in LG(M)$, then $\alpha= d^c \omega^{n-1}$
which gives
$\int_M \alpha \wedge \rho= \int_M \omega^{n-1} \wedge d^c\rho$.
For any non-zero positive current, the integral
$\int_M \omega^{n-1} \wedge d^c\rho$ (known as ``the mass''
of the current) is positive. This means that 
$LG(M)$ is evaluated positively on all elements of ${\cal C}$.
It remains to prove the converse inclusion, that is,
to show that for all $\alpha \in W\backslash LG(M)$,
there exists a closed 1-current $\xi$
such that $d^c \xi$ is positive, but
$\langle \alpha, \xi \rangle =0$.

\hfill

{\bf Step 2:}
Fix $u\in W$ and apply the Hahn--Banach theorem to 
the closed affine space 
\[ V= u + d^c(\text{\sf Aeppli exact $(n-1, n-1)$-forms})+ \text{\sf exact 1-forms}
\]
and the open cone $d^c(\text{\sf strictly positive $(n-1, n-1)$-forms})$.
By Hahn--Banach, these spaces do not intersect if there exists a
$1$-current $\xi$ such that: 
\begin{equation}\label{_3_conditions_HB_LG_Equation_}
	\begin{split}
	\langle -\xi, u + d^c(\text{\sf Aeppli exact $(n-1, n-1)$-forms})\rangle& = 0, \ \ \text{\sf and}\\ 
\langle -\xi, d^c(\text{\sf strictly positive $(n-1, n-1)$-forms})\rangle & >0, \ \ \text{\sf and}\\ 
\langle \xi, \text{\sf exact 1-forms}\rangle & =0. 
	\end{split}
\end{equation}

{\bf Step 3:}
The condition $\langle \xi, \text{\sf exact 1-forms}\rangle =0$ means
that $\xi$ is $d$-closed; indeed, $\langle \xi,
d\zeta\rangle =\pm \langle d\xi,\zeta\rangle.$
However, $d^c(\text{\sf closed 1-forms})\subset \Lambda^{1,1}(M)$,
 because locally a closed 1-form is exact,
and $d^c df\in \Lambda^{1,1}(M)$ for any function $f$.

\hfill

{\bf Step 4:}
Using integration by parts, the second 
condition of \eqref{_3_conditions_HB_LG_Equation_} translates to 
$$\langle d^c\xi, \text{\sf strictly positive $(n-1, n-1)$-forms}\rangle >0,$$
which is equivalent to $d^c\xi$ being a positive current.

\hfill

{\bf Step 5:} The first condition of
\eqref{_3_conditions_HB_LG_Equation_}  implies that
$$\langle -\xi, d^c(\text{\sf Aeppli exact $(n-1, n-1)$-forms})\rangle= 0$$
that is, $\langle d^c\xi, \text{\sf Aeppli exact $(n-1, n-1)$-forms}\rangle =0$.
This is equivalent to $d^c\xi$ being a closed (1,1)-form; it
follows from the same argument that proves Lamari's theorem
(\ref{_Lamari_proved_Theorem_}). 
\endproof

\section{The Lee-Gauduchon cone is bimeromorphically
		invariant}

\definition
	Let $X,Y$ be complex manifolds, and $Z\subset X\times Y$
	a closed subvariety such that the projections of $Z$ to
	$X$ and $Y$ are proper and generically bijective. Then
	$Z$ is called {\bf a bimeromorphic map,} and $X$ and $Y$ are
	called {\bf bimeromorphic.}
	
\hfill

The structure of bimeromorphisms is given in the following fundamental result:

\hfill
	
\theorem {\rm (weak factorization theorem, \cite{_AKMW:weak_factorization_})}\\
	Any bimeromorphism can be decomposed into a composition
		of blow-ups and blow-downs with smooth centers (in arbitrary order).
		
\hfill

\remark This immediately implies that 
	{bimeromorphic manifolds have the same fundamental
		group.} Also, {the spaces of global holomorphic forms
		on bimeromorphic manifolds are naturally isomorphic.}

\hfill

This implies the following

\hfill
	
\corollary
	Let $M_1, M_2$ be compact complex $n$-manifolds which
	are bimeromorphic. {Then $W(M_1)=W(M_2)$,}
	where $W\subset H^{2n-1}(M)$ is the Lee-Gauduchon
        space 
(Subsection \ref{_W:subsection_}).
	
\hfill

\theorem\label{_LG_bimeromo_Theorem_}
	Let $M_1, M_2$ be compact complex $n$-manifolds which
	are bimeromorphic. {Then $LG(M_1)=LG(M_2)$.}\\
	
	\pstep
	Let $[\theta_1], [\theta_2]$ classes in $H^1(M_i)$ 
	which are identified by the natural isomorphism
	$H^1(M_1) = H^1(M_2)$. We are going to prove that
	$[\theta_1]\in {\cal C}(M_1)\Leftrightarrow [\theta_2]\in {\cal C}(M_2)$,
where ${\cal C}(M_i)$ is the dual cone to the LG-cone in $M_i$ (\ref{_LG_dual_to_pos_Theorem_}). 
By definition, ${\cal C}(M_i)$ is the cone of classes  $[\theta_i]\in H^1(M, \R)$ such that
the Bott-Chern class of $d^c([\theta_i])$ is pseudoeffective.
Denote by $\theta_i$ the closed currents representing $[\theta_i]$
such that the currents $d^c \theta_i$ are positive.

\hfill
	
	{\bf Step 2:}
	The pushforward of a positive current is always
	positive; the pullback of a current is, in general, not
	defined. This makes it difficult
	to identify the pseudoeffective cones of
	bimeromorphic manifolds. 
	However, the 1-currents $\theta_i$ are exact on the
	universal cover $\tilde M_i$: $\tilde \theta_i= d f_i$,
	where $f_i$ are generalized functions on $\tilde M_i$.
	If $d^c \theta_i$ is positive, $f_i$ is plurisubharmonic;
	however, a plurisubharmonic function is always
	$L^1_{loc}$-integrable, and can be extended
	over a closed analytic subset,
 \cite[Theorem 5.4, p. 45]{_Demailly:CADG_}. Therefore, the
	pullback {\em and} the pushforward
	of a plurisubharmonic function are plurisubharmonic,
	which implies that $f_1$ can be lifted to the
	graph of the bimeromorphic correspondence
	and pushed forward to a plurisubharmonic
	function on $\tilde M_2$.
	\endproof
	

\section{The Lee--Gauduchon cone: examples}
\label{_examples_Section}


We end this note with results on the Lee--Gauduchon cone on 
special classes of compact complex manifolds: Oeljeklaus-Toma, 
strongly Gauduchon, balanced and LCK.
 
\subsection{Strongly Gauduchon and balanced manifolds}

\definition {(\rm \cite{_Popovici_})} 
A Gauduchon form $\omega^{n-1}$ is called {\bf strongly Gauduchon}
if $\6(\omega^{n-1})$ is $\bar\6$-exact.

\hfill

\remark
This is equivalent to $\omega^{n-1}$ being the 
$(n-1,n-1)$-part of a closed form. Indeed,
$\6(\omega^{n-1})= \bar\6 (x^{n, n-2})$
implies that $d(x^{n, n-2})- d\omega^{n-1}+ d(\overline{x^{n, n-2}})=0$.

\hfill

\remark
The strongly Gauduchon property is implied by
$d\omega^{n-1}=0$; such Hermitian form are known as {\bf
	balanced}, \cite{_Michelson_}.  The class of balanced manifolds
is bimeromorphically invariant, \cite{_Alessandrini_Bassanelli:bimero_}.

\hfill

\proposition
Let $\omega$ be a strongly Gauduchon form. Then
 its Lee-Gauduchon class vanishes.

\hfill

\proof Let $x^{n, n-2}+ \omega^{n-1}+ \overline{x^{n, n-2}}$
be a closed form with $\omega^{n-1}$ a Gauduchon form. Then
$I(x^{n, n-2}+ \omega^{n-1}+ 
\overline{x^{n, n-2}}) =-x^{n, n-2}+ \omega^{n-1} -\overline{x^{n, n-2}}$,
hence $d^c\omega^{n-1} = - d(\omega^{n-1})$.
Therefore, the Lee-Gauduchon class of
	any strongly Gauduchon metric vanishes.
\endproof

\hfill

\remark
If $0\in LG(M)$, then $LG(M)=W$.  Indeed, this cone
is open, convex and $\R^{>0}$-invariant.
Therefore, $LG(M)=W$ for the class of strongly Gauduchon manifolds,
which includes balanced manifolds, which includes
all twistor spaces (of  ASD manifolds,
\cite{_Michelson_,_Hitchin_},  of conformally flat
Riemannian manifolds of even dimension,
\cite{_Gauduchon_twistor_}, of quaternion-K\"ahler
manifolds,  \cite{_Pontecorvo_}, of hyperk\"ahler
manifolds, \cite{_Kaledin_Verbitsky_},  of hypercomplex
manifolds, \cite{_Tomberg_}).

\hfill

\remark
The condition ``$d\omega^{n-1}$ is $d^c$-exact''
is a weaker form of the strongly Gauduchon condition.
A manifold is strongly Gauduchon if and only if
all exact positive (1,1)-currents vanish (\cite{_Popovici_}), and
$d\omega^{n-1}$ is $d^c$-exact if and only if
all positive (1,1)-currents with Bott--Chern 
cohomology classes in $d^c(H^1(M))$
vanish (\ref{_LG_dual_to_pos_Theorem_}).

\hfill

\example
Let $M$ be a Calabi-Eckmann manifold (\cite{_Calabi_Eckmann_}), a complex manifold
diffeomorphic to $S^3\times S^3$ and fibered over $\C P^1 \times \C P^1$
with fiber an elliptic curve. The Bott-Chern group $H^{1,1}_{BC}(M)$
is non-zero; indeed, let $\omega\in \Lambda^{1,1}(M)$ be the pullback of the 
K\"ahler form on $\C P^1 \times \C P^1$. The Bott-Chern class $[\omega]\in H^{1,1}_{BC}(M)$ of this form 
is pseudoeffective, and it is exact, because $H^2(M)=0$. Moreover, $[\omega]\neq 0$,
because a plurisubharmonic function on a compact manifold is constant. Therefore,
$M$ is not strongly Gauduchon; however, the Lee-Gauduchon cone
of $M$ is trivial, because $H^5(M, \R)=0$.

\hfill

\remark
We obtained that for a Calabi-Eckmann manifold,
the cone of positive exact (1,1)-currents and 
the cone of positive (1,1)-currents with Bott--Chern 
cohomology classes in $d^c(H^1(M))$ are distinct.


%
%
%
%
%
%
%
%
%
%
%
%

\subsection{LCK manifolds}

\definition
A Hermitian manifold of complex dimension $ >1$
$(M,I,g,\omega)$ is called 
{\bf locally conformally K\"ahler} (LCK) if 
there exists a closed 1-form $\theta$ such that $d\omega=\theta\wedge\omega$.
The 1-form $\theta$ is called the {\bf Lee form}.

\hfill 

\remark 
We recall that all known compact LCK 
manifolds {belong to one of three classes:
	blow-ups of LCK with potential, blow-ups of 
	Oeljeklaus--Toma and Kato manifolds.} We refer to \cite{_OV_book_} for LCK geometry.

\hfill 

Recently, we proved the following result:

\hfill

\theorem {\rm (\cite{_OV_balanced_})} 
Let $(M, \omega, \theta)$ be a compact LCK manifold in any of
these classes. {Then $d^c \theta$ is exact pseudoeffective.}

	\hfill

\corollary For a compact LCK manifold in any of the above classes, $LG(M)\not\ni 0$ and $M$ does not
admit a strongly Gauduchon metric.

\hfill 

This leads to the following conjecture which was the 
starting point for writing this note:

\hfill

\conjecture
{The two-form $d^c \theta$ is pseudoeffective on all 
	compact LCK manifolds.}
	
\hfill

The Lee-Gauduchon cone for these three
classes of LCK manifolds can be computed 
explicitly.

\hfill

The following proposition was proven in 
\cite{_OV_tohoku_}.

\hfill

\proposition
Let $(M,\omega, \theta)$ be an LCK manifold with
potential (\cite{_OV_book_}). Then 
$H^1(M,\C)= H^{1,0}_d(M) \oplus \overline{H^{1,0}_d(M)}
\oplus \langle\theta\rangle$.

\proof \cite[Theorem 6.1]{_OV_tohoku_}
\endproof

\hfill

From this theorem it is clear that $\dim W=1$, and
since $d^c \theta$ is pseudo-effective on any LCK 
manifold with potential, its Lee-Gauduchon 
cone is a half-line $\{x\in W\ \ |\ \ \langle x, \theta\rangle >0\}$.

The same is true for any Kato manifold $M$
(\cite{_IOP:Kato_,_IOPR:toric_Kato_}).
Indeed, $b_1(M)=1$, because $M$ is a
deformation of a blown-up Hopf manifold,
and $H^{1,0}_d(M)=0$ because 
$H^{1,0}_d(M) \oplus \overline{H^{1,0}_d(M)}$
contribute an even number to $b_1(M)$.
Therefore, $\dim W=1$. Since $d^c \theta$ is pseudo-effective on 
a Kato manifold (\cite{_OV_balanced_}), its Lee-Gauduchon 
cone is a half-line $\{x\in W\ \ |\ \ \langle x, \theta\rangle >0\}$.

We deal with all OT-manifolds, locally conformally K\"ahler
or not, in the next subsection.

\subsection{Oeljeklaus-Toma manifolds}

Oeljeklaus-Toma manifolds were introduced in 
\cite{_Oeljeklaus_Toma_}, and were much studied since then.
These are manifolds associated with number fields.
Their geometry and topology ultimately depends on
two invariants associated with a number field:
its number of real embeddings, denoted $s$,
and the number of complex embeddings, denoted $2t$
(complex embeddings go in pairs, which are complex
conjugate). When $t=1$, an OT manifold is LCK 
(\cite{_Oeljeklaus_Toma_}), and
otherwise it is not (\cite{_Deaconu_Vuletescu_}).
A dimension 2 OT manifold is called {\bf Inoue
surface of type $S^0$}.
OT manifolds are complex solvmanifolds
(\cite{_Kasuya:solvmanifolds_}), that is,
they can be obtained as a left quotient of a solvable
Lie group $G$ with a left-invariant complex structure
by a cocompact lattice $\Gamma$.

One of the main tools which are applied to OT-manifolds
is the averaging procedure, originally defined by F. Belgun
for Inoue surfaces, and then extended to nilmanifolds
in \cite{_Fino_Gra_} (see also \cite{_Kasuya:solvmanifolds_}, \cite{_GFV_algebraic_dimension_}).
Averaging is defined on any quotient manifold
of form $X:=G/\Gamma$, where $G$ is a Lie group and
$\Gamma\subset G$ a discrete subgroup of finite
covolume. The averaging is a map
$\Av:\; \Lambda^k(X) \arrow \Lambda^k(\g^*)$,
where $\Lambda^k(\g^*)$ denotes the space
of antisymmetric $k$-forms on its Lie algebra.
It takes a $k$-form $\alpha$ to a functional 
$T\mapsto \int_{G/\Gamma} \langle T, \alpha\rangle \Vol$,
where $\Vol$ is a left-invariant Haar volume, 
and $T\in \Lambda^k(\g)$ an antisymmetric $k$-vector.
Since $\Lambda^k(\g)^*=\Lambda^k(\g^*)$, the averaging
takes values in $\Lambda^k(\g^*)$.

 The integral $\int_{G/\Gamma} T \wedge \alpha$
is well defined for any $k$-current and any $\alpha\in \Lambda^{n-k}(\g^*)$, hence the averaging is naturally extended to currents
(\cite{_GFV_algebraic_dimension_}).

We will identify $\Lambda^k(\g^*)$ and the space
of left-invariant forms on $G$. The de Rham differential
on $\Lambda^k(\g^*)$ is known as the Chevalley-Eilenberg differential.
It is not hard to see that the averaging commutes with the de
Rham differential; when $G$ is equipped with a left-invariant
complex structure, averaging preserves the Hodge decomposition.
By Nomizu theorem (\cite{_Nomizu:cohomology_}), for  nilmanifolds, the
averaging map induces an isomorphism in cohomology.
This theorem is false for solvmanifolds; however,
H. Kasuya \cite[Lemma 2.1]{_Kasuya:solvmanifolds_} 
has shown that for OT-manifolds 
the averaging induces isomorphism on $H^1$.

The averaging map is very convenient to deal with the 
positive currents: it takes a positive, closed current
to a non-zero, positive, closed invariant form.

The following claim trivially follows from these observations.

\hfill

\claim\label{_averaging_posi_OT_Claim_}
Let $M$ be an OT-manifold, and
$[\alpha]\in H^1(M,\R)$ a cohomology
class; using \cite{_Kasuya:solvmanifolds_},
we represent $[\alpha]$ by a closed, invariant
1-form $\alpha\in \Lambda^1(\g^*)$.
Then $[\alpha]$ can be represented
by a current $\alpha_1$ such that $d^c\alpha_1$ is 
positive if and only if $d^c\alpha$ is positive.

\hfill

\proof
The form $\Av(\alpha_1)$ is equal to $\alpha$,
because these forms are cohomologous, and 
$d:\;\Lambda^0(\g^*)\arrow  \Lambda^1(\g^*)$ is zero.
Since the averaging preserves positivity,
the (1,1)-form $d^c \alpha= \Av(d^c(\alpha_1))$ is also positive.
\endproof

\hfill

To proceed, we need to describe the
complex structure on the Lie algebra of the
solvable group $G$ associated with the OT-manifold
explicitly.

An OT-manifold is constructed from a number
field $K$ which has $2t$ complex embeddings and 
$s$ real ones (in other words, $K\otimes_\Q \R= \C^{2t}\oplus \R^s$).
For such a number field, the group of units $\calo_K^*$ has rank
$t+s-1$ by Dirichlet unit theorem. Oeljeklaus and
Toma construct a torsion-free subgroup $U\subset \calo_K^*$ 
of rank $s$, such that the quotient described below
is compact.

Let $\tilde M_0= \C^t\otimes \R^s$; we identify $\tilde M$
with a subspace of $K\otimes_\Q \R= \C^{2t}\oplus \R^s$
obtained by taking only one $\C$ in each pair of
complex conjugate components. The additive group 
$\calo_K= \Z^{2s+t}$ acts on $\tilde M_0$ by translations;
this action is cocompact, totally discontinuous
and the quotient $\frac{\tilde M_0}{\calo_K}$
is a torus.

Now, let ${\Bbb H}$ denote the upper half-plane.
Consider the manifold $\tilde M:= \C^t\otimes {\Bbb H}^s$,
identified with $\tilde M_0\times \R^s$.
We equip  the quotient manifold  $\frac {\tilde M}{\calo_K}$
with the action of $U$
as follows. The action of $U$ on $\tilde M_0= K\otimes_\Q \R$
is the standard multiplicative action; it clearly
commutes with the $\calo_K$-action on the same manifold.
To describe the action of $U$ on each ${\Bbb H}$-component,
we index these components by the set 
$\{\sigma_1, ... \sigma_s\}$ of real
embeddings of $K$. Then $u\in U$ acts on 
the ${\Bbb H}$-component number $k$
as a multiplication by the real number $\sigma_k(u)$.

It is not hard to see that if we choose
generators $u_i$ of $U$ in such a way that
$\sigma_k(u_i) >1$ for each $k, i$, then the
quotient $M =\frac {\tilde M}{\calo_K\rtimes U}$
is compact; the action of $\calo_K\rtimes U$
on $\C^t\otimes {\Bbb H}^s\subset \C^{s+t}$ 
is by construction complex affine. Clearly, $M$
is a cocompact quotient of the Lie group
$G:=(\calo_K\otimes_\Z\R)\rtimes (U\otimes_\Z\R)$
obtained as a semidirect product of the 
abelian Lie groups identified with the vector spaces
$\calo_K\otimes_\Z\R$ and $U\otimes_\Z\R$.

Following \cite[Section 6]{_Kasuya:solvmanifolds_},
we fix the basis of the dual space to the Lie algebra $\g$ of $G$.
Let $\gamma_i$ be generators associated with the $\C$-factors,
with the complex structure $I$ taking
$\gamma_{2i-1}$ to $\gamma_{2i}$,
and $\alpha_i, \beta_i$ be the generators
associated with the ${\Bbb H}$-factors, with 
$I(\alpha_i)=\beta_i$. The standard coordinate
system on ${\Bbb H}$ is $x=\Re(z), y=\log(\Im(z))$;
then $\beta_i$ correspond to the real components,
$\beta_i=\frac{dx}y$ and $\alpha_i$ correspond to the imaginary
component, $\alpha_i=\frac{dy}{y}$.

In this basis, the Chevalley-Eilenberg differential (the operator
dual to the Lie bracket) is written as in \cite[Section 6]{_Kasuya:solvmanifolds_}:
$d\alpha_i=0, d\beta_i=-\alpha_i\wedge \beta_i$, and
\[ d\gamma_{2i-1}= \psi_i \wedge \gamma_{2i-1} +
\phi_i \wedge \gamma_{2i}, \ \ d\gamma_{2i}= - \phi_i \wedge \gamma_{2i-1} +
\psi_i \wedge \gamma_{2i},
\]
where $\psi_i, \phi_i$ are linear combinations of $\alpha_i$.

Since $I$ exchanges $\alpha_i$ with $\beta_i$,
the intersection $d\g^* \cap I(d\g^*)$ is generated
by $d \beta_i$, with $d\gamma_i$ contributing nothing.
However, $-d\beta_i$ are positive for all $i$.
This brings the following proposition.

\hfill

\proposition\label{_pos_exact_OT_Proposition_}
Let $M= \frac{\C^t\times {\Bbb H}^s}{\calo_K\rtimes U}$ be 
an OT-manifold, $z_i$ the complex coordinates on the ${\Bbb H}$-components,
and $d^c d(\log(\Im z_i))$ the K\"ahler form associated
with the Poincar\'e metric on this ${\Bbb H}$-component.
Then the cone $A$ of pseudo-effective Bott--Chern cohomology classes 
which belong to  $d^c(H^1(M,\R))$ is generated by positive linear
combinations of $d^c d(\log(\Im z_i))$.

\hfill

\proof
By \ref{_averaging_posi_OT_Claim_},
the cone $A$ is generated by positive forms $d^c\rho$,
where $\rho\in \Lambda^1(\g^*)$ is closed.
However, the intersection $d\Lambda^1(\g^*)\cap I(d \Lambda^1(\g^*))$
is generated by $\alpha_i\wedge \beta_i$, 
as indicated above, hence the only positive
(1,1)-forms in $d\Lambda^1(\g^*)$ are 
linear combinations of 
$-d\beta_i$.

By definition, $\beta_i$ is the $G$-invariant form on
the $i$-th ${\Bbb H}$-component of $\C^t\times {\Bbb H}^s$,
corresponding to $d^c\log \Im(z_i)$, hence $d^c d(\log(\Im z_i))=-d\beta_i$.
\endproof

\hfill

Using \ref{_pos_exact_OT_Proposition_} in 
conjunction to \ref{_LG_dual_to_pos_Theorem_}, we can compute the Lee-Gauduchon
cone of an OT-manifold $M$ explicitly. Recall that $H^1(M)$
is generated by $d\log\Im z_i$, where $z_i$ are
coordinates on the ${\Bbb H}$-components of $\C^t\times {\Bbb H}^s$
(\cite{_Oeljeklaus_Toma_,_Kasuya:solvmanifolds_}).
However, the forms $d^c \alpha_i$ are positive.
By \ref{_LG_dual_to_pos_Theorem_}, the
Lee-Gauduchon cone of an OT-manifold is dual
to the cone of cohomology classes $\alpha$ in
$H^1(M)$ such that $d^c\alpha$ is positive.
This brings the following result.

\hfill

\theorem
Let $M$ be an OT-manifold, 
and $[\alpha_1], ..., [\alpha_s]$ generators
of the first cohomology constructed above,
$\alpha_i=d\log\Im z_i$. Then 
the LG-cone of $M$ is dual to the convex
cone generated by $\alpha_i$ in $H^1(M,\R)$.
\endproof



\hfill

\noindent{\bf Acknowledgements:}
We are grateful to Dan Popovici for an important 
email consversation and consultations about the 
strong Gauduchon metrics.

\hfill

\hfill

{\scriptsize

	\noindent {\sc Liviu Ornea\\
		{\sc University of Bucharest, Faculty of Mathematics and Informatics, \\14
			Academiei str., 70109 Bucharest, Romania}, \\
		also:\\
		Institute of Mathematics ``Simion Stoilow" of the Romanian
		Academy\\
		21, Calea Grivitei Str.
		010702-Bucharest, Romania}\\
	{\tt lornea@fmi.unibuc.ro,   liviu.ornea@imar.ro}
	
	\hfill

	\noindent
	{\sc Misha Verbitsky\\
		{\sc Instituto Nacional de Matem\'atica Pura e
			Aplicada (IMPA) \\ Estrada Dona Castorina, 110\\
			Jardim Bot\^anico, CEP 22460-320\\
			Rio de Janeiro, RJ - Brasil }\\
		{\tt verbit@impa.br}
}}
\end{document}